\documentclass{amsart}
\usepackage{graphicx}
\title{Generalised Mycielski graphs and bounds on chromatic numbers}
\author[G. Simons, C. Tardif, D. Wehlau]{Gord Simons, Claude Tardif, 
David Wehlau}\thanks{the third author's research is supported by grants 
from NSERC and ARP}
\address{Royal Military College of Canada \\ PO Box 17000 Stn Forces,
Kingston, ON\\ Canada, K7K 7B4}

\newcommand{\cqfd}{\begin{flushright}\rule{8pt}{9pt}\end{flushright} \par}
\newtheorem{define}{Definition}
\newtheorem{proposition}[define]{Proposition}

\newtheorem{problem}[define]{Problem}
\newtheorem{theorem}[define]{Theorem}
\newtheorem{lemma}[define]{Lemma}
\newtheorem{remark}[define]{Remark}
\newtheorem{corollary}[define]{Corollary}

\newcommand{\bd}{\begin{define} \rm}
\newcommand{\ed}{\end{define}}

\newcommand{\coind}{\mathrm{coind}}

\newcommand{\boco}[1]{\mathrm{B}(#1)}

\newfont{\Bb}{msbm10 scaled\magstep1}

\begin{document}

\begin{abstract}
We prove that the coindex of the box complex $\boco{H}$ of a graph $H$ can be measured
by the generalised Mycielski graphs which admit a homomorphism to it.
As a consequence, we exhibit for every graph $H$ a system of linear equations
solvable in polynomial time, with the following properties:
If the system has no solutions, then $\coind(\boco{H}) + 2 \leq 3$;
if the system has solutions, then $\chi(H) \geq 4$. We generalise
the method to other bounds on chromatic numbers using linear algebra.
\end{abstract} 

\maketitle
{\small \noindent{\bf Keywords:}
Graph colourings, homomorphisms, Box complexes, generalised Mycielski graphs
\newline \noindent {\em AMS 2010 Subject Classification:}
05C15.}
\smallskip

\section{Introduction}
For any integer $k \geq 2$ and real number $\epsilon \in (0,1)$, the {\em Borsuk graph} 
$B_{k,\epsilon}$ is the graph whose vertices are the points of the $(k-2)$-sphere
$S_{k-2} \subseteq \mathbb{R}^{k-1}$, and whose edges join pairs of points 
$X, Y$ that are ``almost antipodal''
in the sense that the norm of $X-Y$ is at least $2 - \epsilon$. 
In \cite{eh}, Erd\H{o}s and Hajnal
used the Borsuk-Ulam theorem to prove that the chromatic
number of $B_{k,\epsilon}$ is $k$.
In fact, they proved that the statement $\chi(B_{k,\epsilon}) = k$ is equivalent
to the Borsuk-Ulam theorem. For some years this result remained a curiosity
involving infinite graphs. Then Lov\'asz~\cite{lovasz} devised complexes
that allow the use of the Borsuk-Ulam Theorem to find lower bounds on 
chromatic numbers of finite graphs, and used this method to prove the Kneser
conjecture on the chromatic number of the Kneser graphs.

Lov\'asz' method inspired many adaptations and developments, giving rise to the
field of ``topological lower bounds'' on the chromatic number of a graph.
Our work is inspired by a bound in terms of ``coindices of box complexes'', 
specifically
$$\chi(H) \geq \coind(\boco{H}) + 2.$$
The relevant definitions of the coindex $\coind(\boco{H})$ of the box complex
$\boco{H}$ of $H$ are well detailed in \cite{matousek,mz,st}. However, our intent 
is to avoid the topological setting. We will use the following result of Simonyi and Tardos.
\begin{theorem}[\cite{st}] \label{sita}
For any graph $H$, $\coind(\boco{H}) + 2$ is the largest $k$ such that there exist
$\epsilon > 0$ for which $B_{k,\epsilon}$ admits a homomorphism (that is,
an edge-preserving map) to $H$. 
\end{theorem}
This result indeed allows us to restrict our discussion to 
the field of graphs and homomorphisms: We can alternatively define $\coind(\boco{H}) + 2$
{\em as} the largest $k$ such that there exist an $\epsilon > 0$ for which $B_{k,\epsilon}$ 
admits a homomorphism to $H$. 
This viewpoint yields an economy of definitions, but is not necessarily practical
for computational purposes. Indeed in many cases a knowledge of simplicial complexes
and topological tricks is needed to compute $\coind(\boco{H}) + 2$
and effectively bound $\chi(H)$.

Dochtermann and Schultz \cite{ds} found finite (``spherical'') graphs that 
play the role of the Borsuk graphs in Theorem~\ref{sita}. 
In this note we show that the generalised Mycielski graphs can also be used in the
the same role. We do not know whether this alternative presentation yields
effective computations in the general case. However, for low values of $\coind(\boco{H})$,
our definition indeed leads to practical calculations that can be shown to be
conclusive in some cases. The method, in turn, inspires effectively computable 
lower bounds on the chromatic number of a graph.

\section{Generalised Mycielski graphs}

We will use the following definitions of categorical products, 
looped paths, cones and generalised Mycielski graphs.
The {\em categorical product} of two graphs $G$ and $G'$ is the graph $G\times G'$ defined by
\begin{eqnarray*}
V(G\times G') & = & V(G) \times V(G'), \\
E(G \times G') & = & \{ [(u,u'),(v,v')] : [u,v]  \in E(G) \mbox{ and } [u',v'] \in E(G') \}.
\end{eqnarray*}
We sometimes use directed graphs as factors, and view undirected 
graphs as symmetric directed graphs.
In this case, all square brackets representing edges in the above 
definition should be replaced by parentheses representing arcs.
For $n \in \mbox{\Bb N}^*$, the {\em looped path} $\mbox{\Bb P}_n$ is the path 
with vertices $0, 1, \ldots, n$ linked consecutively, with a loop at $0$. 
For a graph $G$, the {\em $n$-th cone} (or {\em $n$-th generalised Mycielskian})
$M_n(G)$ over $G$ is the graph 
$$(G \times \mbox{\Bb P}_n)/ \sim_n,$$
where $\sim_n$ is the equivalence which identifies all vertices whose second coordinate is $n$.
The classes $\mathcal{K}_k$ of {\em generalised Mycielski graphs} are defined
recursively as follows: $\mathcal{K}_2 = \{ K_2 \}$, and for $k \geq 3$,
$$
\mathcal{K}_k = \{ M_n(G) : G \in \mathcal{K}_{k-1}, n \in \mbox{\Bb N}^*\}.
$$

Csorba \cite{csorba1, csorba2} proved that for any graph $H$ and integer $n$,
the geometric realisation of $\boco{M_n(H)}$ is $\mbox{\Bb Z}_2$-homotopy
equivalent to the geometric realisation of the suspension of $\boco{H}$. 
In particular this implies that for every $G \in \mathcal{K}_k$
we have
$$\coind(\boco{G}) + 2 = \chi(G) = k.$$ 
(Recall that $\coind(\boco{G})$ is defined in \cite{matousek,mz,st} and characterised
implicitly by Theorem~\ref{sita} above.)

\begin{lemma} \label{gminb}
For every Borsuk graph $B_{k,\epsilon}$, there exists a graph $G$
in $\mathcal{K}_k$ such that $G$ admits a homomorphism to $B_{k,\epsilon}$.
\end{lemma}

\smallskip \noindent {\em Proof.}
$\mathcal{K}_2 = \{K_2\}$, and $B_{2,\epsilon} = K_2$
for every $\epsilon > 0$. Suppose that $G \in \mathcal{K}_{k-1}$ admits a homomorphism
to $B_{k-1,\epsilon/2}$. Put $n = \lceil \pi/\epsilon \rceil$.
We will show that $M_n(G)$ admits a homomorphism to $B_{k,\epsilon}$
We identify
the vertex set of $B_{k-1,\epsilon/2}$ with the equator of $B_{k,\epsilon}$.
Hence $[u,v] \in E(B_{k-1,\epsilon/2})$ implies $[u,v] \in E(B_{k,\epsilon})$.
Let $p_N, p_S$ respectively be the north and south poles of $B_{k,\epsilon}$.
Let $\phi: G \rightarrow B_{k-1,\epsilon/2}$ be a homomorphism. 
For every $u \in V(G)$, let $\phi(u) = u_{N,0}, u_{N,1}, \ldots, u_{N,n} = p_N$
be equally spaced points on the quarter of the great circle joining $\phi(u)$ and $p_N$.
Similarly let  $\phi(u) = u_{S,0}, u_{S,1}, \ldots, u_{S,n} = P_S$
be equally spaced points on the quarter of the great circle joining $\phi(u)$ and $P_S$.
Define $\psi: M_m(G) \rightarrow  B_{k,\epsilon}$ by
$$
\psi(u,i) = \left \{
\begin{array}{l}
\mbox{$u_{N,i}$ if $i$ is even,}\\
\mbox{$u_{S,i}$ if $i$ is odd.}
\end{array} 
\right.
$$
Note that $\psi(u,n) = p_N$ or $p_S$ according to whether $n$ is even or odd, hence
$\psi$ is well defined. Also, $\psi$ extends $\phi$. For $[u,v] \in E(G)$ and $i<m$, we have 
\begin{eqnarray*}
||\psi(u,i) + \psi(v,i+1)|| & = & ||u_{N,i} + v_{S,i+1}|| = ||u_{N,0} + v_{S,1}|| \\
& \geq & ||u_{N,0} + v_{S,0}|| - ||v_{S,0} - v_{S,1}|| > 2 - \epsilon.
\end{eqnarray*}
Therefore $\psi$ is a homomorphism. \cqfd

\begin{corollary} \label{gmgtc}
For any graph $H$, $\coind(\boco{H}) + 2$ is the largest $k$ such that there exists
a $G$ in $\mathcal{K}_k$ admitting a homomorphism to $H$. 
\end{corollary}

\smallskip \noindent {\em Proof.} Let $k = \coind(\boco{H}) + 2$.
By Theorem~\ref{sita} and Lemma~\ref{gminb}, there exist
a number $\epsilon > 0$ and a graph $G \in \mathcal{K}_k$
such that there are homomorphisms of $B_{k,\epsilon}$ to $H$
and of $G$ to $B_{k,\epsilon}$. The composition of these is
a homomorphism of $G$ to $H$. On the other hand, it is well known
that if $G$ admits a homomorphism to $H$, then
$\coind(\boco{G}) \leq \coind(\boco{H})$.
\cqfd

The usefulness of the bound $\chi(H) \geq \coind(\boco{H}) + 2$ derives
from the fact that the chromatic number is hard to compute, hence
lower bounds are useful. However for a finite graph $H$, $\chi(H)$
can at least be determined by a finite computation, while 
$\coind(\boco{H})$ is not known to be computable.

The class $\mathcal{K}_3$ consists of the odd cycles. 
Therefore the problem of determining whether an input graph $H$ satisfies 
$\coind(\boco{H}) + 2 \leq 2$ is equivalent to that of determining 
whether $H$ is bipartite, which admits an efficient solution. 
In the remainder of the paper, we focus on the implication
$$\chi(H) \leq 3 \Rightarrow \coind(\boco{H}) + 2 \leq 3.$$
We present an approach derived from Corollary~\ref{gmgtc}.

\section{Signatures of odd cycles}\label{sigoc}

We begin by reinterpreting homomorphisms of cones in terms of 
paths in exponential graphs. We first present the basic properties
of exponential graphs which can be found in standard references,
e.g.~\cite{hn}. For two graphs $G$ and $H$, the {\em exponential graph} 
$H^G$ has for vertices all functions $f: V(G) \rightarrow V(H)$, 
and for edges all pairs $[f,g]$ of functions
such that for every $[u,v] \in E(G)$, $[f(u),g(v)] \in E(H)$. In particular 
$f$ is a homomorphism of $G$ to $H$ if and only if $f$ is a loop in $H^G$.
A homomorphism of $G \times G'$ to $H$ corresponds to a homomorphism
of $G'$ to $H^G$. In particular, with $G' = \mbox{\Bb P}_m$, we have the following.
\begin{remark} \label{egc}
A homomorphism of $M_m(G)$ to $H$ corresponds to an $m$-path in
$H^G$ from a loop to a constant map.
\end{remark}
We will suppose that $H$ is connected, thus all the constant maps
are in the same connected component of $H^G$, which we call the
{\em connected component of the constants}.
Suppose that $H$ satisfies $\coind(\boco{H})+2 \geq 4$. 
Then by Corollary~\ref{gmgtc}, there exists
an odd cycle $C$ such that for some integer $m$, $M_m(C) \in \mathcal{K}_4$
admits a homomorphism to $H$. By the above remark, this is equivalent to
the existence of a loop in the connected component of the constants
in $H^C$. Since $H^C$ is finite, such a loop could be found in finite time. 
Thus for a fixed $C$, the question of the existence of a $m$ such that 
$M_m(C)$ admits a homomorphism to $H$ is decidable. 
However to settle the question
as to whether $\coind(\boco{H})+2 \geq 4$, we need to look at the infinitely
many possible choices for $C$, so the search is still infinite.

Our next step is to look for easily computable
invariants that are constant on the components of $H^C$.
For a graph $H$, let $A(H)$ denote the set of its arcs. That is,
for $[u,v] \in E(H)$, $A(H)$ contains the two arcs $(u,v)$ and $(v,u)$.
Let  $\mathcal{F}^*(A(H))$ be the free group generated 
by the elements of $A(H)$. 
Let $C$ be an odd cycle with vertex-set $\mathbb{Z}_{2n+1} = \{0, \ldots, 2n\}$.
The elements of $H^C$ are not necessarily homomorphisms of $C$ to $H$.
However, if $(f,g)$ is an arc of $H^C$, then $f(0), g(1), f(2), \ldots, g(2n)$
is a closed walk of length $4n+2$ in $H$. We define the 
{\em $\mathcal{F}^*(A(H))$-signature}
$\sigma_{\mathcal{F}^*(A(H))}(f,g)$ of $(f,g)$ as
$$
\sigma_{\mathcal{F}^*(A(H))}(f,g) 
= \Pi_{i = 0}^{2n} (f(2i),g(2i+1)) \cdot (f(2i+2),g(2i+1))^{-1}.
$$
The indices are of course taken modulo $2n+1$. Also since the product
is noncommutative, it is necessary to specify that the product is developed 
left to right: $\Pi_{i = 0}^{k}x_i = x_0 x_1 \cdots x_k$ 
rather than $x_k x_{k-1} \cdots x_0$.
Thus $\sigma_{\mathcal{F}^*(A(H))}(f,g)$ is essentially a list of the arcs
in the closed walk generated by $f$ and $g$, with trivial simplifications.

We define a congruence $\theta^*$ on $\mathcal{F}^*(A(H))$ as follows:
If $u, v, w, x$ is a $4$-cycle of $H$, we put  
$$
\begin{array}{rcl}
(u,v) \cdot (w,v)^{-1} & \theta^* & (u,x) \cdot (w,x)^{-1},
\end{array}
$$
that is,
$$
\begin{array}{rcl}
(u,v) \cdot (w,v)^{-1} \cdot (w,x) \cdot (u,x)^{-1} & \theta^* & 1_{\mathcal{F}(A(H))}.
\end{array}
$$
The group $\mathcal{G}^*(H)$ is defined by
$$\mathcal{G}^*(H) := \mathcal{F}^*(A(H))/ \theta^*.$$

Now if $(f,g)$ and $(f,g')$ are two arcs of $H^C$, then the 
corresponding terms in the product defining their signatures are
congruent:
$$
(f(2i),g(2i+1)) \cdot (f(2i+2),g(2i+1))^{-1}
\theta^*
(f(2i),g'(2i+1)) \cdot (f(2i+2),g'(2i+1))^{-1}
$$
for $i = 0, \ldots, 2n$.
Therefore $\sigma_{\mathcal{F}^*(A(H))}(f,g)/ \theta^*
= \sigma_{\mathcal{F}^*(A(H))}(f,g')/ \theta^*$.
We define the {\em $\mathcal{G}^*(H)$-signature} $\sigma_{\mathcal{G}^*(H)}(f)$
of a non-isolated vertex $f$ in $H^C$ by
$$\sigma_{\mathcal{G}^*(H)}(f) = \sigma_{\mathcal{F}^*(A(H))}(f,g)/ \theta^*$$
for any neighbour $g$ of $f$.

\begin{proposition} \label{signature}
Let $C$ be an odd cycle and $f: C \rightarrow H$ a homomorphism.
If $f$ belongs to the connected component of the constants in $H^C$, then
$$\sigma_{\mathcal{G}^*(H)}(f) = 1_{\mathcal{G}^*(H)}.$$
\end{proposition}

\smallskip \noindent {\em Proof.}
If $f$ is a constant map in $H^C$, then for any neighbour $g$ of $f$ we have
\begin{eqnarray*}
\sigma_{\mathcal{F}^*(A(H))}(f,g) 
& = & \Pi_{i = 0}^{2n} (f(2i),g(2i+1)) \cdot (f(2i+2),g(2i+1))^{-1} \\
& = & \Pi_{i = 0}^{2n} (f(0),g(2i+1)) \cdot (f(0),g(2i+1))^{-1}
= 1_{\mathcal{F}^*(A(H))},
\end{eqnarray*}
hence $\sigma_{\mathcal{G}^*(H)}(f) = 1_{\mathcal{G}^*(H)}$. 

To extend the argument by connectivity,
we will use a natural automorphism of $\mathcal{G}^*(H)$.
Define $\alpha_0^*: A(H) \rightarrow \mathcal{F}^*(A(H))$ by 
$\alpha_0^*(u,v) = (v,u)^{-1}$.
Then $\alpha_0^*$ extends to an order 2 automorphism $\alpha^*$ of $\mathcal{F}^*(A(H))$.
For every arc $(f,g)$ of $H^C$, the product defining $\sigma_{\mathcal{F}^*(A(H))}(g,f)$
uses the term $\alpha^*(e)$ for every term $e$ used in the product defining
$\sigma_{\mathcal{F}^*(A(H))}(f,g)$. More precisely 
$\sigma_{\mathcal{F}^*(A(H))}(g,f) = \alpha^*(W') \alpha^*(W)$, where
$W$ is the product of the $2n+1$ first terms in $\sigma_{\mathcal{F}^*(A(H))}(f,g)$
and $W'$ is the product of the $2n+1$ last terms. In particular,
$\sigma_{\mathcal{F}^*(A(H))}(g,f)$ is a conjugate of 
$\alpha^*(\sigma_{\mathcal{F}^*(A(H))}(f,g)) = \alpha^*(W) \alpha^*(W')$. 
Now for a generator
$y = (u,v) \cdot (w,v)^{-1} \cdot (w,x) \cdot (u,x)^{-1}$ of 
$1_{\mathcal{F}^*(A(H))}/ \theta^*$, $\alpha(y)$ is a conjugate of the generator
$(x,u) \cdot (v,u)^{-1} \cdot (w,v) \cdot (w,x)^{-1}$ of
$1_{\mathcal{F}^*(A(H))}/ \theta^*$. Hence $1_{\mathcal{F}^*(A(H))}/ \theta^*$
is invariant under $\alpha^*$, whence $\alpha^*$ induces an automorphism of
$\mathcal{G}^*(H)$, which we also call $\alpha^*$.

Thus if $f$ is any element of $H^C$ such that
$\sigma_{\mathcal{G}^*(H)}(f) = 1_{\mathcal{G}^*(H)}$, then for any neighbour $g$ of $f$,
$\sigma_{\mathcal{G}^*(H)}(g)= \sigma_{\mathcal{F}^*(A(H))}(g,f)/ \theta^*$
is a conjugate of 
$\alpha^*(\sigma_{\mathcal{F}^*(A(H))}(f,g)/ \theta^*) = 1_{\mathcal{G}^*(H)}$.
Therefore $\sigma_{\mathcal{G}^*(H)}(g) = 1_{\mathcal{G}^*(H)}$.
By connectivity, this implies that $\sigma_{\mathcal{G}^*(H)}(f)$
is identically $1_{\mathcal{G}^*(H)}$ on the connected component of the
constants in $H^C$. \cqfd

\begin{figure}[h]
\centering
\includegraphics{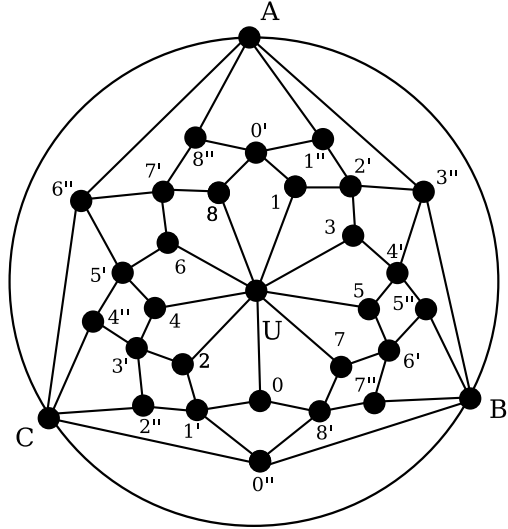}
\caption{Witness to the fallacy of the converse of Proposition~\ref{signature}.}  
\label{example}
\end{figure}

The following example shows that the converse of Proposition~\ref{signature}
does not hold in general. Consider the graph $H$ in Figure~\ref{example}. 
Let $f: C_3 \rightarrow H$
be the homomorphism defined by $f(0) = A$, $f(1) = C$ and $f(2) = B$.
Any neighbour $g$ of $f$ is in the set $S$ defined by
$$
S = \{ g \in H^{K_3} : g(0) \in \{ A, 0''\}, g(1) \in \{C, 3''\}, g(2) \in \{B, 6''\} \}.
$$
It is easy to check that any neighbour of an element of $S$ is again in $S$.
Thus $f$ is not in the connected component of the constants in $H^{C_3}$.
We will show that $\sigma_{\mathcal{G}^*(H)}(f) = 1_{\mathcal{G}^*(H)}$.

Let $C_9$ be the $9$-cycle with vertex-set $\mbox{\Bb Z}_9$.
Let $h_0, h_1, h_2, h_3, h_4: \mbox{\Bb Z}_9 \rightarrow V(H)$ be defined
as follows: $h_0(i) = U$ for all $i \in \mbox{\Bb Z}_9$, $h_1, h_2, h_3$ are defined by 
$h_1(i) = i$, $h_2(i) = i'$, $h_3(i) = i''$, and $h_4$ is defined by
$$
(h_4(0), h_4(1), \ldots, h_4(8)) = (A,C,A,C,B,C,B,A,B).
$$
Since $h_0$ is a constant and $h_0, h_1, h_2, h_3, h_4$ is a path in $H^{C_9}$,
we have $\sigma_{\mathcal{G}^*(H)}(h_4) = 1_{\mathcal{G}^*(H)}$. Note that
$h_4$ is a homomorphism, hence 
$\sigma_{\mathcal{G}^*(H)}(h_4) = \sigma_{\mathcal{F}^*(A(H))}(h_4,h_4)/\theta^*$.
By definition we have
\begin{eqnarray*}
\sigma_{\mathcal{F}^*(A(H))}(h_4,h_4) & = &
(A,C) \cdot (A,C)^{-1} \cdot (A,C) \cdot (B,C)^{-1} \cdot (B,C) \cdot (B,C)^{-1} \\
& & \cdot (B,A) \cdot (B,A)^{-1} \cdot (B,A) \cdot (C,A)^{-1} \cdot (C,A) \cdot (C,A)^{-1} \\
& & \cdot (C,B) \cdot (C,B)^{-1} \cdot (C,B) \cdot (A,B)^{-1} \cdot (A,B) \cdot (A,B)^{-1} \\
& = & (A,C) \cdot (B,C)^{-1} \cdot (B,A) \cdot (C,A)^{-1} \cdot (C,B) \cdot (A,B)^{-1} \\
& = & \sigma_{\mathcal{F}^*(A(H))}(f,f).
\end{eqnarray*}
Therefore $\sigma_{\mathcal{G}^*(H)}(f) = \sigma_{\mathcal{F}^*(A(H))}(f,f)/\theta^* 
= 1_{\mathcal{G}^*(H)}$.

Note that $h_4: C_9 \rightarrow H$ factors as $f \circ f'$, with $f': C_9 \rightarrow C_3$ given
by
$$
(f'(0),f'(1), \ldots, f'(8)) = (0,1,0,1,2,1,2,0,2).
$$
Therefore $h_4$ could be seen as an ``unfolding'' of $f$, which falls
in the connected component of the constants. It is not clear whether a similar
phenomenon always occurs.
\begin{problem} \label{naconverse}
Let $H$ be a graph, $C_n$ an odd cycle and $f$ a homomorphism in $H^{C_n}$
such that $\sigma_{\mathcal{G}^*(H)}(f) = 1_{\mathcal{G}^*(H)}$.
Does there exist an odd cycle $C_m$ and a homomorphism $f': C_m \rightarrow C_n$
such that $f \circ f'$ is in the connected component of a constant in $H^{C_m}$?
\end{problem}
If Problem~\ref{naconverse} has an affirmative answer, then detecting the existence of
a homomorphism of some $M_m(C') \in \mathcal{K}_4$ which admits a homomorphism
is equivalent to finding an odd closed walk in $H$ with trivial signature.
We do not know of a feasible approach to the latter problem. However we
will see that the Abelian relaxation of the problem is tractable.

Let $\gamma$ be the commutator of $\mathcal{G}^*(H)$.
The group $\mathcal{G}(H)$ is defined as $\mathcal{G}^*(H)/ \gamma$,
and the {\em abelian signature} of $f \in H^C$ is defined
as $\sigma_{\mathcal{G}(H)}(f) \equiv \mathcal{G}^*(H)/ \gamma$.
Thus $\mathcal{G}(H) = \mathbb{Z}^{A(H)}/ \theta$, where $\theta$ is generated
by the relations
$$
\begin{array}{rcl}
(u,v) - (w,v) + (w,x) - (u,x) & \theta & 0_{\mathbb{Z}^{A(H)}}.
\end{array}
$$
such that $u, v, w, x$ is a $4$-cycle of $H$, and 
$$
\sigma_{\mathcal{G}(H)}(f) = 
\left ( \sum_{i = 0}^{2n} (f(2i),g(2i+1)) - (f(2i+2),g(2i+1)) \right ) / \theta,
$$
where $g$ is any neighbour of $f$.
As a consequence of Proposition~\ref{signature}, we have the following
\begin{corollary} \label{absignature}
Let $C$ be an odd cycle and $f: C \rightarrow H$ a homomorphism.
If $f$ belongs to the connected component of the constants in $H^C$, then
$$\sigma_{\mathcal{G}(H)}(f) = 0_{\mathcal{G}(H)}.$$ \cqfd
\end{corollary}
In the next section we see that the search for $C$ and $f$ that satisfy the 
conclusion of Corollary~\ref{absignature} is tractable. 

\section{The signature system of equations}\label{sigsys}

Let $H$ be a connected graph. To each arc $(u,v)$ of $H$ we associate
an integer variable $X_{u,v}$. The {\em flow constraint} at a vertex
$u$ of $H$ is the equation 
\begin{equation}\label{flow}
\sum_{v \in N_H(u)} (X_{u,v} - X_{v,u}) = 0.
\end{equation}
We also consider a {\em parity constraint} requiring that the sum of these
variables is odd:
\begin{equation}\label{parity}
\sum_{(u,v) \in A(H)} X_{u,v} - 2N = 1.
\end{equation}
Thus we introduce an additional integer variable $N$. Finally, the 
{\em signature constraint} is the equation
\begin{equation}\label{nullsign}
\left ( \sum_{(u,v) \in A(H)} (X_{u,v} - X_{v,u})\cdot (u,v) \right ) / \theta 
 = 0_{\mathcal{G}(H)}.
\end{equation}

\begin{figure}[htp]
\centering
\includegraphics{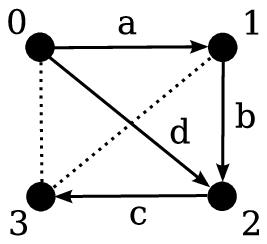}
\caption{$K_4$}  
\label{k4}
\end{figure}

Note that while the flow and parity constraints are equations in $\mathbb{Z}$,
the signature constraint is an equation in $\mathcal{G}(H)$. However it can 
be rewritten as a system of equations representing the coordinate
of the vectors in the finitely generated abelian group $\mathcal{G}(H)$. 
For instance consider $K_4$ with vertex-set $\{ 0, 1, 2, 3\}$ as illustrated in
figure~\ref{k4}. 

The arcs $(0,1), (1,2), (2,3)$ and $(0,2)$ are denoted $a, b, c, d$ respectively.
For $e = (i,j)$, we will write $e^{-}$ for $(j,i)$.  
By definition of $\theta$ we have
\begin{eqnarray*}
(0,3) / \theta & = & (a - b^{-} + c) / \theta, \\
(3,0) / \theta & = & (c^{-} - b + a^{-}) / \theta, \\
(1,3)/ \theta & = & (b - d + (0,3)) / \theta = (b - d + a - b^{-} + c) / \theta,\\ 
(3,1)/ \theta & = & (c^{-} - d + a) / \theta, \\
(2,0)/ \theta & = & (c - (1,3) + a^{-}) / \theta 
= (b^{-} - a + d -b + a^{-}) / \theta.
\end{eqnarray*}
It is easy to check that for $B = \{a, a^{-}, b, b^{-}, c, c^{-}, d\}$,
the natural homomorphism of $\mathbb{Z}^B$ to $\mathcal{G}(K_4)$ is an isomorphism,
thus we can identify $\mathcal{G}(K_4)$ with $\mathbb{Z}^B$.
Therefore the signature constraint on $K_4$ can be rewritten as
$$
z_a \cdot a + z_{a^{-}} \cdot a^{-}
+ z_b \cdot b + z_{b^{-}} \cdot z^{-}
+ z_c \cdot c + z_{c^{-}} \cdot c^{-}
+ z_d \cdot d = 0,
$$
which only has the trivial solution. In terms of the variables
$X_{i,j}$, this yields the following system of equations.
\begin{eqnarray*}
z_a = 0 & : & (X_{0,1}-X_{1,0}) + (X_{0,3}-X_{3,0}) - (X_{2,0}-X_{0,2})= 0, \\
z_{a^{-}} = 0 & : & (X_{1,0}-X_{0,1}) + (X_{3,0}-X_{0,3}) + (X_{2,0}-X_{0,2})= 0, \\
 z_b = 0 & : & (X_{1,2}-X_{2,1}) - (X_{3,0}-X_{0,3}) 
+ (X_{1,3}-X_{3,1}) - (X_{2,0}-X_{0,2})= 0, \\
z_{b^{-}} = 0 & : & (X_{2,1}-X_{1,2}) - (X_{0,3}-X_{3,0}) - 
(X_{1,3}-X_{3,1}) + (X_{2,0}-X_{0,2})= 0, \\
z_c = 0 & : & (X_{2,3}-X_{3,2}) + (X_{0,3}-X_{3,0}) + (X_{1,3}-X_{3,1}) = 0, \\
 z_{c^{-}} = 0 & : & (X_{3,2}-X_{2,3}) + (X_{3,0}-X_{0,3}) + (X_{3,1}-X_{1,3}) = 0, \\
 d = 0 & : &  0 = 0.
\end{eqnarray*}

The {\em signature system} of $H$ is the system of linear equations consisting 
of the flow constraint~(\ref{flow}) at every vertex of $H$, the parity 
constraint~(\ref{parity}) and the signature constraint~(\ref{nullsign}).
\begin{proposition} \label{linsyst}
Let $H$ be a connected graph. 
Then the signature system of $H$ admits integers solutions 
if and only if there exists an odd cycle
$C$ and a homomorphism $f: C \rightarrow H$ such that 
$\sigma_{\mathcal{G}(H)}(f) = 0_{\mathcal{G}(H)}$.
\end{proposition}

\smallskip \noindent {\em Proof.}
Let $C$ be an odd cycle (with $V(C) = \mathbb{Z}_{2n+1}$) and 
$f: C \rightarrow H$ a homomorphism. Put $N = n$ and
$$X_{u,v} = |\{ i \in \mathbb{Z}_{2n+1} : f(i) = u, f(i+1) = v \}|$$
for every arc $(u,v) \in A(H)$.
Then the set of values $X_{u,v}, (u,v) \in A(H)$ and 
$N$ satisfy the flow constraints~(\ref{flow})
and the parity constraint~(\ref{parity}). We have
$$\sigma_{\mathcal{G}^(H)}(f) = \left (\sum_{(u,v) \in A(H)} (X_{u,v} - X_{v,u})\cdot (u,v) \right )/ \theta,$$
and this value is $0_{\mathcal{G}(H)}$ if and only if the signature
constraint is satisfied, that is, $X_{u,v}, (u,v) \in A(H)$ and $N$
are solutions to the signature system of $H$.

Conversely, let $X_{u,v}, (u,v) \in A(H)$ and $N$ be an integer solution to the system.
We first modify the solution by subtracting $\min\{X_{u,v}, X_{v,u}\}$ from both
$X_{u,v}$ and $X_{v,u}$ for every edge $[u,v]$ of $H$ and subtracting the sum of these minima
from $N$. We now have a non-negative integer solution.
We then add $1$ to every variable $X_{u,v}, (u,v) \in A(H)$, and $|E(H)|$ to $N$.
This yields a positive integer solution with connected support.
Let $G$ be the multidigraph with $V(G) = V(H)$ and $X_{u,v}$ parallel arcs connecting
$u$ to $v$ for every $(u,v) \in A(H)$. Then $G$ is connected Eulerian, and an Euler closed
trail in $G$ corresponds to a homomorphism $f: C \rightarrow H$ with $|V(C)| = 2N+1$
and 
$$\sigma_{\mathcal{G}^(H)}(f) = 
\left (\sum_{(u,v) \in A(H)} (X_{u,v} - X_{v,u})\cdot (u,v) \right ) / \theta 
= 0_{\mathcal{G}(H)}.$$
\cqfd

For instance, let us return to the example of $K_4$ discussed above.
It is easy to see that every homomorphism $f$ of an odd cycle $C$ to $K_4$
is at distance at most two from any constant in $K_4^C$, hence 
$\sigma_{\mathcal{G}(K_4)}(f) = 0_{\mathcal{G}(K_4)}$. Therefore the signature
system of $K_4$ should not be that hard to solve. In fact we see that 
the condition $z_d = 0$ of the signature constraint is trivially satisfied.
We recognize the flow constraint at vertex $3$ in the condition $z_{c^{-\alpha}} = 0$,
and moreover all conditions of the signature constraint
reduce to flow constraints with elementary manipulations. 
Thus the signature constraint is redundant on $K_4$,
and every solution to the flow and parity constraints 
is a solution to the signature system.

By Corollary~\ref{gmgtc}, Remark~\ref{egc} and Corollary~\ref{absignature}, 
Proposition~\ref{linsyst} has the following consequence.
\begin{corollary} \label{coss}
If the signature system of $H$ has no integer solutions,
then 
$$\coind(\boco{H}) + 2 \leq 3.$$ 
\cqfd
\end{corollary}

\begin{figure}[h]
\centering
\includegraphics{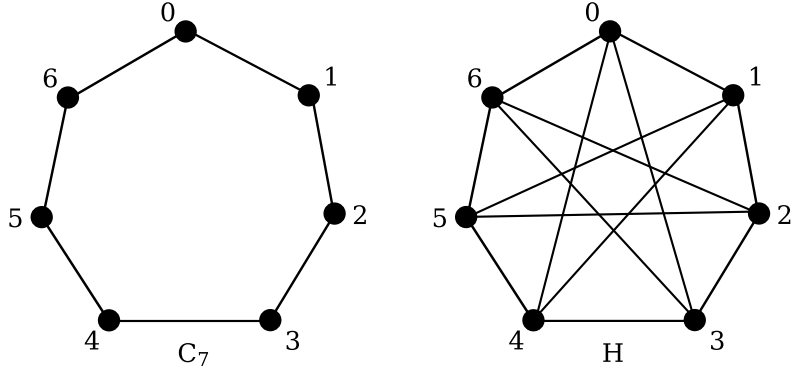}
\caption{Two applications of Corollary~\ref{coss}}  
\label{c7andh}
\end{figure}

For instance, the seven-cycle $C_7$ with vertex-set $\mathbb{Z}_7$ shown in 
Figure~\ref{c7andh} has no $4$-cycles, thus 
$\mathcal{G}(C_7) = \mathbb{Z}^{A(C_7)}$.
Therefore the signature constraint on $C_7$ implies that
$X_{i,i+1} = X_{i+1,i}$ for all $i \in \mathbb{Z}_7$, 
which is incompatible with the parity constraint.
By Corollary~\ref{coss}, this implies that 
$\coind(\boco{C_7}) + 2 \leq 3$, which is a well-known fact.

The second graph $H$ in Figure ~\ref{c7andh} is the ``third relational
power'' of $C_7$ obtained by adding an edge between the extreme points 
of the 3-paths. Thus we have
\begin{eqnarray*}
(i,i+3) / \theta & = & (i,i+1) - (i+2,i+1) + (i+2,i+3) / \theta, \\
(i+3,i) / \theta & = & (i+3,i+2) - (i+1,i+2) + (i+1,i) / \theta
\end{eqnarray*}
for all $i \in \mathbb{Z}_7$. It is easy to check that the relations
derived from the $4$-cycles of the form $i, i+1, i+4, i+3$ are redundant.
Therefore we again have $\mathcal{G}(H) = \mathbb{Z}^{A(C_7)}$.
Hence the signature constraint on $H$ is of the form
$$
\sum_{i\in \mathbb{Z}_7} z_{i,i+1} \cdot (i,i+1) 
+  \sum_{i\in \mathbb{Z}_7} z_{i+1,i} \cdot (i+1,i) 
 =  0.
$$
From it we get the constraints
\begin{eqnarray*}
z_{i,i+1} = 0 & : & (X_{i,i+3}-X_{i+3,i}) + (X_{i-1,i+2}-X_{i+2,i-1}) \\
&  & + (X_{i-2,i+1}-X_{i+1,i-2}) + (X_{i,i+1}-X_{i+1,i}) = 0.
\end{eqnarray*}
Adding these for $i = 0, \ldots, 6$ we get
$$
\sum_{i\in \mathbb{Z}_7} \left ( 3(X_{i,i+3}-X_{i+3,i}) + (X_{i,i+1}-X_{i+1,i}) \right ) =  0.
$$
We then add the parity constraint to get
$$
\sum_{i\in \mathbb{Z}_7} \left ( 4X_{i,i+3}- 2X_{i+3,i} + 2X_{i,i+1} \right ) - 2N =  1,
$$
which has no integer solutions. Therefore the signature system of $H$ is inconsistent,
and by Corollary~\ref{coss}, this implies that 
$\coind(\boco{H}) + 2 \leq 3$. (This is again a well-known fact.)

Our next example shows that the converse of Corollary~\ref{coss} does not hold.
We consider the graph $U(5,3)$ studied in~\cite{stv, zimmerman}. The vertices
of $U(5,3)$ are the ordered pairs $(i,\{j,k\})$ such that $i, j, k \in \{1, \ldots, 5\}$
and $i \not \in \{j, k\}$. Two of these vertices $(i,\{j,k\})$, $(i',\{j',k'\})$
are joined by an edge if $i \in \{j',k'\}$ and $i' \in \{j,k\}$. 
In~\cite{stv}, Simonyi, Tardos and Vre\'{c}ica proved that 
$\coind(\boco{U(5,3)}) + 2 = 3$. Using a different topological bound,
they also proved that $\chi(U(5,3)) = 4$.

$U(5,3)$ has 90 edges thus 180 arcs, and 105 $4$-cycles.
With such larger graphs, it is useful to use a computer algebra
package to ease computations. In~\cite{zimmerman}, Zimmerman used 
the package Magma~\cite{magma} to prove that $\mathcal{G}(U(5,3))$ 
is isomorphic to $\mathbb{Z}^{71}$, and that the signature system of $U(5,3)$
admits integer solutions. In particular this shows that the
converse of Corollary~\ref{coss} does not hold.

The signature system of a graph $H$ is an efficiently solvable
system that gives information about $\coind(\boco{H})$
when the system has no solution. Our next result shows that
signature system still gives information about $H$ even when it 
has a solution.

\begin{proposition} \label{4chrom}
If the signature system of $H$ has integer solutions,
then $\chi(H) \geq 4$. 
\end{proposition}

We will prove Proposition~\ref{4chrom} as a consequence of
Proposition~\ref{homsign} of the next section. To summarize,
the problem of determining whether a graph $H$ satisfies
$\coind(\boco{H}) + 2 \geq 3$ inspired the method of the signature
system. In turn this method is an efficiently computable criterion
to prove that $\chi(H) \geq 4$, that includes and expands the cases
where $\coind(\boco{H}) + 2 \geq 3$. Also, this method can be generalised,
as shown in the next section.
  
\section{Generalised signature systems}

A {\em valued digraph} $(D,\phi)$ is a directed graph $D$
along with an integer valuation $\phi: A(D) \rightarrow \mathbb{Z}$
of its arcs. For a graph $H$, the {\em congruence generated by $(D,\phi)$} 
is the congruence $\theta(D,\phi)$ on $\mathbb{Z}^{A(H)}$ 
generated by the conditions
$$
\left ( \sum_{(u,v) \in A(D)} \phi(u,v) \cdot (f(u),f(v)) \right ) / \theta(D,\phi)
 =  0_{\mathbb{Z}^{A(H)}} / \theta(D,\phi)
$$
for every homomorphism $f$ of $D$ to $H$. 

\begin{figure}[h]
\centering
\includegraphics{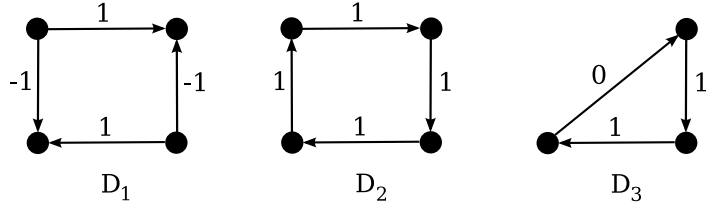}
\caption{Valued digraphs}  
\label{valdi}
\end{figure}

In Figure~\ref{valdi}, $(D_1,\phi_1)$ models our original congruence:
$\theta = \theta(D_1,\phi_1)$. Note that homomorphisms of $D_1$ to $H$ may
collapse vertices. However this only gives the trivial condition
$0 = 0$ in $\theta(D_1,\phi_1)$. However, collapsing $D_2$ on an edge $[u,v]$
of $H$ gives the condition 
$$
2((u,v) + (v,u)) / \theta(D_2,\phi_2) = 0_{\mathbb{Z}^{A(H)}} / \theta(D_2,\phi_2).
$$
Thus $\mathbb{Z}^{A(H)}/ \theta(D_2,\phi_2)$ always has torsion.
We have not encountered an example where 
$\mathcal{G}(H) = \mathbb{Z}^{A(H)}/ \theta(D_1,\phi_1)$ has torsion,
though we cannot prove that $\mathcal{G}(H)$ is always torsion-free.

Collapsing $D_2$ on two incident edges $[u,v]$, $[v,w]$
of $H$ gives the condition
$$
(u,v) + (v,u) \theta(D_2,\phi_2) (v,w) + (w,v).
$$
By connectivity, the terms $(u,v) + (v,u)$ are then congruent for all edges 
$[u,v]$ of $H$ (or on the edges of a connected component if $H$ is not connected).
Therefore on a genuine $4$-cycle $u, v, w, x$ of $H$, we have
$$
(u,v) + (v,w) + (w,x) + (x,u)/ \theta(D_2,\phi_2) = 0_{\mathbb{Z}^{A(H)}} / \theta(D_2,\phi_2)
$$
and 
$$
(v,w) + (w,v) + (x,u) + (u,x)/ \theta(D_2,\phi_2) = 0_{\mathbb{Z}^{A(H)}} / \theta(D_2,\phi_2).
$$
Subtracting the second relation from the first, we get
$$
(u,v) - (w,v) + (w,x) - (u,x) \theta(D_2,\phi_2) = 0_{\mathbb{Z}^{A(H)}} / \theta(D_2,\phi_2).
$$
Thus on any graph $H$ we have $\theta(D_1,\phi_1) \subseteq \theta(D_2,\phi_2)$.

For a set $\mathcal{D}$ of valued digraphs, we put 
$\theta(\mathcal{D}) = \vee_{(D,\phi) \in \mathcal{D}} \theta(D,\phi)$.
For a graph $H$, the group $\mathcal{G}_{\mathcal{D}}(H)$ is defined
as $\mathbb{Z}^{A(H)}/ \theta(\mathcal{D})$.
The $\mathcal{D}$-signature constraint on integer variables $X_{u,v}, (u,v) \in A(H)$
is the condition
\begin{equation}\label{nullsignd}
\left ( \sum_{(u,v) \in A(H)} (X_{u,v} - X_{v,u})\cdot (u,v) \right ) / \theta(\mathcal{D}) 
 = 0_{\mathcal{G}_{\mathcal{D}}(H)},
\end{equation}
and the $\mathcal{D}$-signature system is the system consisting of the 
flow constraints~(\ref{flow}), the parity constraint~(\ref{parity}) 
and the $\mathcal{D}$-signature constraint~(\ref{nullsignd}).

\begin{proposition} \label{homsign}
Let $H$, $H'$ be  graphs such that there 
exists a homomorphism $\psi: H \rightarrow H'$.
If the $\mathcal{D}$-signature system has integer solutions on $H$, 
then it has integer solutions on $H'$. 
\end{proposition}

\smallskip \noindent {\em Proof.} 
First note that $\hat{\psi}: A(H) \rightarrow A(H')$ defined
by $\hat{\psi}(u,v) = (\psi(u),\psi(v))$ extends to a group homomorphism 
of $\mathbb{Z}^{A(H)}$ to $\mathbb{Z}^{A(H')}$, which we also denote 
$\hat{\psi}$. We will show that $\hat{\psi}$ maps 
$0_{\mathbb{Z}^{A(H)}} / \theta(\mathcal{D})$
to $0_{\mathbb{Z}^{A(H')}} / \theta(\mathcal{D})$.
Indeed any generator of $0_{\mathbb{Z}^{A(H)}} / \theta(\mathcal{D})$
is of the form $\sum_{(u,v) \in A(D)} \phi(u,v) \cdot (f(u),f(v))$,
where $(D,\phi)$ is in $\mathcal{D}$ and $f$ is a homomorphism of $D$ in $H$.
It is mapped by $\psi$ to 
$\sum_{(u,v) \in A(D)} \phi(u,v) \cdot (\psi \circ f(u),\psi \circ f(v))$.
The latter is in $0_{\mathbb{Z}^{A(H')}} / \theta(\mathcal{D})$,
since $\psi \circ f$ is a homomorphism of $D$ to $H'$. Therefore
$\hat{\psi}$ induces a group homomorphism of 
$\mathbb{Z}^{A(H)} / \theta(\mathcal{D}) = \mathcal{G}_{\mathcal{D}}(H)$ to 
$\mathbb{Z}^{A(H')} / \theta(\mathcal{D}) = \mathcal{G}_{\mathcal{D}(H')}$,
which we again call $\hat{\psi}$.

Now suppose that $X_{u,v}, (u,v) \in A(H)$ and $N$ are solutions to the 
$\mathcal{D}$-signature system of $H$. For $(u',v') \in A(H')$, put
$$
X'_{u',v'} = \sum \{ X_{u,v} : \psi(u,v) = u' \mbox{ and } \psi(v) = v'\}.
$$

\smallskip \noindent {\em Claim 1.} $X'_{u',v'}, (u',v') \in A(H')$ and $N$
satisfy the flow constraint~(\ref{flow}) on $H'$.

Indeed for $u' \in V(H')$ we have
$$
\sum_{v' \in N_{H'}(u')} (X'_{u',v'} - X'_{v',u'}) 
= \sum_{u \in \psi^{-1}(u')} \sum_{v \in N_H(u)} (X_{u,v} - X_{v,u}) = 0.
$$

\smallskip \noindent {\em Claim 2.} $X'_{u',v'}, (u',v') \in A(H')$ and $N$
satisfy the parity constraint~(\ref{parity}) on $H'$.

Indeed we have
$$
\sum_{(u',v') \in A(H')} X'_{u',v'} - 2N
= \sum_{(u,v) \in A(H)} X_{u,v} - 2N = 1.
$$

\smallskip \noindent {\em Claim 3.} $X'_{u',v'}, (u',v') \in A(H')$ and $N$
satisfy the $\mathcal{D}$-signature constraint~(\ref{nullsignd}) on $H'$.

Indeed since
$$
\left ( \sum_{(u,v) \in A(H)} (X_{u,v} - X_{v,u})\cdot (u,v) \right ) / \theta(\mathcal{D})
=  0_{\mathcal{G}_{\mathcal{D}}(H)}
$$
and $\hat{\psi}$ is a group homomorphism, we have
$$
\left ( \sum_{(u,v) \in A(H)} (X_{u,v} - X_{v,u})\cdot (\psi(u),\psi(v)) \right ) 
/ \theta(\mathcal{D})
=  0_{\mathcal{G}_{\mathcal{D}}(H')}.
$$
Regrouping preimages we get
$$
\left ( \sum_{(u',v') \in A(H')} (X'_{u',v'} - X'_{v',u'})\cdot (u',v') \right ) 
/ \theta(\mathcal{D})
=  0_{\mathcal{G}_{\mathcal{D}}(H')}.
$$

Claims 1, 2 and 3 prove that
$X'_{u',v'}, (u',v') \in A(H')$ and $N$ are integer solutions 
to the $\mathcal{D}$-signature system on $H'$. \cqfd

The usefulness of Proposition~\ref{homsign} resides in the cases where it can be shown
that the $\mathcal{D}$-signature system has no solutions on a fixed target graph
$H'$. The $\mathcal{D}$-signature system then provides a criterion for the existence
 of a homomorphism of an input graph $H$ to $H'$. In the case of 
Proposition~\ref{4chrom} of the previous section, we have $H' = K_3$
and $\theta(\mathcal{D}) = \theta$.

\smallskip \noindent {\em Proof of Proposition~\ref{4chrom}.}

Since $K_3$ has no $4$-cycles, $\mathcal{G}(K_3) = \mathbb{Z}^{A(K_3)}$. 
Therefore the signature constraint~(\ref{nullsign})
$$
\left ( \sum_{(u',v') \in A(K_3)} (X'_{u',v'} - X'_{v',u'})\cdot (u',v') \right ) / \theta 
=  0_{\mathcal{G}(K_3)}
$$
implies $X'_{u',v'} = X'_{v',u'}$ for all $(u',v') \in A(K_3)$. We then have
$\sum_{(u',v') \in A(H')} X'_{u',v'}$ even, which is incompatible with the
parity constraint~(\ref{parity}). Therefore the signature system has no 
integer solutions on $K_3$. Therefore by Proposition~\ref{homsign}, if the
signature system has integer solutions on a graph $H$,
then $H$ admits no homomorphism to $K_3$, whence $\chi(H) \geq 4$. \cqfd

The signature system was used to prove that the second graph $H$ of 
Figure~\ref{c7andh} satisfies $\coind(\boco{H}) + 2 \leq 3$.
We now use $\mathcal{D} = \{ (D_2,\phi_2), (D_3,\phi_3) \}$
to prove that $\chi(H) \geq 4$.

Let $e_1, e_2, e_3$ be the three clockwise arcs in a triangle of $H$.
The conguence generated by $(D_3,\phi_3)$ implies
$$
e_1 / \theta(\mathcal{D}) = - e_2 / \theta(\mathcal{D}) =
e_3 / \theta(\mathcal{D}) = -e_1 / \theta(\mathcal{D}).
$$
Thus $e_1, e_2$ and $e_3$ are all congruent to the
same element of order $2$ in $\mathcal{G}_{\mathcal{D}}(H)$.
Extending the argument to all triangles in $H$, we get that
all the arcs $(i,i+1)$ and $(i,i+3)$ are congruent to the
same element of order $2$  in $\mathcal{G}_{\mathcal{D}}(H)$
and similarly all the arcs $(i+1,i)$ and $(i+3,i)$ are congruent to the
same element of order $2$  in $\mathcal{G}_{\mathcal{D}}(H)$.
Now since $\mathcal{D}$ contains $(D_2,\phi_2)$, we get
$$
\left ( (0,1) + (1,2) + (2,3) + (3,0) \right ) / \theta(\mathcal{D})
= 0_{\mathcal{G}_{\mathcal{D}}(H)}.
$$
Therefore $\mathcal{G}_{\mathcal{D}}(H) \simeq \mathbb{Z}_2$, with
all arcs congruent to the non-zero element.
The $\mathcal{D}$-signature constraint on $H$ then reduces to the 
trivial condition $0 = 0$. Any odd cycle corresponds to a solution 
to the $\mathcal{D}$-signature system on $H$, so this system admits 
non-trivial solutions on $H$.

However, $\mathcal{G}_{\mathcal{D}}(K_3) \simeq \mathbb{Z}_2^2$,
with $(0,1), (1,2), (2,0)$ congruent to a non-zero element $a$ of
$\mathcal{G}_{\mathcal{D}}(K_3)$ and $(0,2), (2,1), (1,0)$ congruent to 
another non-zero element $b$ of $\mathcal{G}_{\mathcal{D}}(K_3)$.
The $\mathcal{D}$-signature constraint on $K_3$ is
\begin{eqnarray*}
\left ( \sum_{i \in \mathbb{Z}_3} (X_{i,i+1} - X_{i+1,i}) \right ) \cdot a
+  \left ( \sum_{i \in \mathbb{Z}_3} (X_{i+1,i} - X_{i,i+1}) \right )\cdot b
& = & 0_{\mathcal{G}_{\mathcal{D}}(H)}.
\end{eqnarray*}
It is satisfied only when $ \sum_{i \in \mathbb{Z}_3} (X_{i,i+1} - X_{i+1,i})$
is even, and this is incompatible with the parity constraint. Thus the
$\mathcal{D}$-signature system has no solution on $K_3$. Therefore
$H$ admits no homomorphism to $K_3$, hence $\chi(H) \geq 4$.

We note that generalised signature systems always win the day,
albeit in a trivial way:
\begin{remark} \label{wtd}
Let $H, H'$ be graphs such that $H$ is not bipartite and
there is no homomorphism of $H$ to $H'$. Then there exists a set $\mathcal{D}$
of valued digraphs such that the  $\mathcal{D}$-signature system admits solutions
on $H$ but not on $H'$.
\end{remark}

\smallskip \noindent {\em Proof.} Put $\mathcal{D} = \{ (H,\phi) \}$, where
$\phi$ has value $1$ on the forward arcs of an odd cycle $C$ of $H$, $-1$ on the
backward arcs of $C$ and $0$ elsewhere. Then $C$ corresponds to 
a solution to the $\mathcal{D}$-signature system on $H$, while 
$\mathcal{G}_{\mathcal{D}}(H') = \mathbb{Z}^{A(H')}$ so that the
$\mathcal{D}$-signature system has no solution on $H'$. \cqfd

However, generalised signature systems can be said to be efficient only
with ${\mathcal{D}}$ fixed and $H$ variable. However the model
could deviate even further from the topological problem that inspired
it: signature constraint could be modified
to any condition of the form
\begin{equation}\label{nullsigndmod}
\left ( \sum_{(u,v) \in A(H)} p X_{u,v} + q X_{v,u})\cdot (u,v) \right ) / \theta(\mathcal{D}) 
 = 0_{\mathcal{G}_{\mathcal{D}}(H)}
\end{equation}
with $p, q \in \mathbb{Z}$, and Claim 3 of the proof of Proposition~\ref{homsign}
would remain valid. 
Likewise the parity constraint could be replaced by any condition of the form
\begin{equation}\label{paritymod}
\sum_{(u,v) \in A(H)} X_{u,v} - pN = q,
\end{equation}
and Claim 2 of the proof of Proposition~\ref{homsign} would remain valid. 
Consider a system with an arbitrary set of such constraints, perhaps 
involving many different sets $\mathcal{D}_i$ of valued graphs and 
many variables for each arc. If such a system is solvable on $H$ and
$H$ admits a homomorphism to $H'$, then the system is solvable 
on $H'$.

Nonetheless, with all these generalisations available, there remains
work to be done. It is easy to see that the signature system admits
solutions on any $4$-chromatic generalised Mycielski graph.
Therefore by Proposition~\ref{4chrom}, we have $\chi(G) \geq 4$
for all $G \in \mathcal{K}_4$. By Lemma~\ref{gminb}
this implies $\chi(B_{k,\epsilon}) \geq 4$ for all $\epsilon > 0$.
This is a proof of the Borsuk-Ulam Theorem for the $2$-sphere.
It would be interesting to know whether generalised signature
systems can be used to prove all of the Borsuk-Ulam theorem
by graph-theoretic methods. Perhaps the next step would be to settle
the following. 
\begin{problem}
Does there exist a set $\mathcal{D}$ of valued digraphs such that
if the $\mathcal{D}$-signature system is solvable on a graph $H$,
then $\chi(H) \geq 5$, and otherwise $\coind(\boco{H}) + 2\leq 4$?
\end{problem}


\begin{thebibliography}{99}
\bibitem{magma} W. Bosma, J. Cannon, and C. Playoust, 
The Magma algebra system. I. The user language, 
J. Symbolic Comput. 24(3-4) (1997), 235--265. Computational algebra
and number theory (London, 1993).
\bibitem{csorba1} P. Csorba, 
Fold and Mycielskian on homomorphism complexes,
Contrib. Discrete Math. 3 (2008), 1--8. 
\bibitem{csorba2} P. Csorba,
On the simple $\mbox{\Bb Z}_2$-homotopy types of graph complexes 
and their simple $\mbox{\Bb Z}_2$-universality,
Canad. Math. Bull. 51 (2008), 535--544. 
\bibitem{ds} A. Dochtermann, C. Schultz, Topology of Hom complexes and test graphs
for bounding chromatic number, Israel J. Math. 187 (2012), 371--417.
\bibitem{eh} P. Erd\H{o}s, A. Hajnal, On chromatic graphs, Mat. Lapok 18 (1967), 1--4.
\bibitem{hn} P. Hell, J. Ne\v{s}et\v{r}il
Graphs and homomorphisms, Oxford Lecture Series in Mathematics
 and Its Applications, vol.~28, Oxford University Press, 2004. 
\bibitem{lovasz} L. Lov\'{a}sz, Kneser's conjecture, chromatic number, and homotopy,
J. Combin. Theory Ser. A 25 (1978), 319--324.
\bibitem{matousek} J. Matou\v{s}ek,
{\em Using the Borsuk-Ulam theorem},
Lectures on topological methods in combinatorics and geometry, 
Written in cooperation with Anders Bj\"{o}rner and G\"{u}nter M. Ziegler, 
Universitext, Springer-Verlag, Berlin, 2003. xii+196 pp.
\bibitem{mz} J, Matou\v{s}ek, G. Ziegler, 
Topological lower bounds for the chromatic number: a hierarchy, 
Jahresber. Deutsch. Math.-Verein. 106 (2004), 71--90. 
\bibitem{st} G. Simonyi, G. Tardos, Local chromatic number, Ky Fan's theorem, and circular colorings, Combinatorica 26 (2006), 587--626.
\bibitem{stv} G. Simonyi, G. Tardos, S. Vre\'{c}ica, 
Local chromatic number and distinguishing the strength of topological obstructions,
Trans. Amer. Math. Soc. 361 (2009),  889--908. 
\bibitem{zimmerman} S. Zimmerman, Topological and Algebraic Lower Bounds on
the Chromatic Number of Graphs, M. Sc. Thesis, Queen's University, 2014.

\end{thebibliography}
\end{document}